%&amstex     
\input amstex\documentstyle {amsppt}  
\pagewidth{12.5 cm}\pageheight{19 cm}\magnification\magstep1
\topmatter
\title Character sheaves on disconnected groups, III\endtitle
\author G. Lusztig\endauthor
\address Department of Mathematics, M.I.T., Cambridge, MA 02139\endaddress
\thanks Supported in part by the National Science Foundation\endthanks
\endtopmatter   
\document
\define\mpb{\medpagebreak}
\define\co{\eta}
\define\hco{\hat\eta}
\define\Up{\Upsilon}
\define\Upo{\Upsilon^\o}
\define\baSi{\bar\Si}
\define\haG{\hat\G}
\define\bag{\bar\g}
\define\da{\dagger}

\define\uce{\un{\ce}}

\define\Lie{\text{\rm Lie }}

\define\frl{\forall}

\define\si{\sim}
\define\wt{\widetilde}
\define\sqc{\sqcup}

\define\qua{\quad}

\define\hG{\hat G}

\define\hM{\hat M}

\define\bY{\bar Y}

\define\op{\oplus}

\define\em{\emptyset}
\define\imp{\implies}

\define\iy{\infty}
\define\m{\mapsto}
\define\do{\dots}

\define\lra{\leftrightarrow}

\define\sub{\subset}

\define\T{\times}
\define\ti{\tilde}
\define\nl{\newline}
\redefine\i{^{-1}}
\define\fra{\frac}
\define\un{\underline}

\define\ot{\otimes}
\define\bbq{\bar{\QQ}_l}

\define\tr{\text{\rm tr}}

\define\supp{\text{\rm supp}}

\define\a{\alpha}
\redefine\b{\beta}
\redefine\c{\chi}
\define\g{\gamma}
\redefine\d{\delta}
\define\e{\epsilon}

\redefine\o{\omega}
\define\p{\pi}
\define\ph{\phi}
\define\ps{\psi}

\define\s{\sigma}
\redefine\t{\tau}

\define\k{\kappa}
\redefine\l{\lambda}
\define\z{\zeta}

\redefine\G{\Gamma}

\define\Si{\Sigma}

\define\boc{\bold c}

\define\kk{\bold k}

\redefine\AA{\bold A}

\define\FF{\bold F}

\define\QQ{\bold Q}

\define\cd{\Cal D}
\define\ce{\Cal E}
\define\cf{\Cal F}

\define\ch{\Cal H}

\define\cl{\Cal L}

\define\car{\Cal R}
\define\cs{\Cal S}

\define\cu{\Cal U}
\define\cv{\Cal V}
\define\cw{\Cal W}
\define\cz{\Cal Z}

\define\fc{\frak c}

\define\fg{\frak g}

\define\ft{\frak t}

\define\fK{\frak K}

\define\tY{\ti Y}
\define\tZ{\ti Z}

\define\bS{\bar S}

\define\bce{\bar\ce}
\define\tce{\ti\ce}
\define\CSII{L8}
\define\ADI{L9}
\define\ADII{L10}

\head Introduction\endhead
Throughout this paper, $G$ denotes a fixed, not necessarily connected, reductive 
algebraic group over an algebraically closed field $\kk$. This paper is a part of a
series (beginning with \cite{\ADI},\cite{\ADII}) which attempts to develop a theory
of character sheaves on $G$.

Assume now that $\kk$ is an algebraic closure of a finite field $\FF_q$ and that $G$
has a fixed $\FF_q$-structure with Frobenius map $F$. Let $(L,S,\ce,\ph_0)$ be a
quadruple where $L$ is an $F$-stable Levi of some parabolic of $G^0$, $\ce$ is a 
local system on an isolated $F$-stable stratum $S$ of $N_GL$ with certain properties
and $\ph_0$ is an isomorphism of $\ce$ with its inverse image under the Frobenius 
map. To $(L,S,\ce)$ we have associated in 5.6 an intersection cohomology complex 
$\fK=IC(\bY,\p_!\tce)$ on $G$. Moreover, $\ph_0$ gives rise to an isomorphism $\ph$
between $\fK$ and its inverse image under the Frobenius map. There is an associated
characteristic function $\c_{\fK,\ph}$ (see 15.12(a)) which is a function 
$G^F@>>>\bbq$, constant on $(G^0)^F$-conjugacy classes. The main result of this 
paper is Theorem 16.14 which shows that the computation of this function can be 
reduced to an analogous computation involving only unipotent elements in a smaller 
group (the centralizer of a semisimple element). (This is a generalization of 
\cite{\CSII, Theorem 8.5}. However, even if $G$ is assumed to be connected, as in
\cite{\CSII}, our Theorem 16.14 is more general than that in \cite{\CSII}, since
here we do not make the assumption that $\ce$ is cuspidal. Also, unlike the proof in
\cite{\CSII}, the present proof does not rely on the classification of cuspidal 
local systems.) A main ingredient in Theorem 16.14 are the generalized Green 
functions, see 15.12(c), which generalize those in \cite{\CSII, 8.3.1}. One of the
key properties of the generalized Green functions is the invariance property
15.12(d). In the connected case, this property was stated in \cite{\CSII, 8.3.2},
but the proof given there was incomplete (as pointed out by F.Letellier). Most of
Section 15 is devoted to establishing this invariance property. 

\head Contents\endhead
15. Generalized Green functions.

16. The characteristic function $\c_{\fK,\ph}$.

\head 15. Generalized Green functions\endhead
\subhead 15.1 \endsubhead
In this section we fix a pair $(L,\Si)$ where $L$ is a subgroup of $G^0$ and $\Si$
is a non-empty subset of $N_GL$; we assume that there exists a parabolic of $G^0$
normalized by $\Si$, with Levi $L$, and that $\Si=\cup_{j\in J}S_j$ where $S_j$ are
distinct isolated strata of $N_GL$ with $\dim S_j$ independent of $j$. Several 
definitions in Sections 3 and 5 which concern the special case where $\Si$ is a 
single stratum will now be extended to the present, more general case. Let 
$\baSi,\bS_j$ be the closure of $\Si,S_j$ in $N_GL$. Then 
$\baSi=\cup_{j\in J}\bS_j$. Since $\bS_j$ is a union of isolated strata of $G$, the
same holds for $\baSi$. Let $\ce$ be a local system on $\Si$ such that 
$\ce_j:=\ce|_{S_j}\in\cs(S_j)$ for all $j\in J$. Since $\Si$ is a smooth, open, 
dense subvariety of pure dimension of $\baSi$, $IC(\baSi,\ce)\in\cd(\baSi)$ is well
defined. Let
$$\Si^*=\{g\in\Si;Z_G(g_s)^0\sub L\}=\cup_{j\in J}S_j^*$$
($S_j^*=\{g\in S_j;Z_G(g_s)^0\sub L\}$ as in 3.11). Then $\Si^*$ is an open dense 
subset of $\Si$ (this follows from the fact that $S_j^*$ is open dense in $S_j$ for
any $j$, see 3.11). Let
$$\tY_{L,\Si}=\{(g,xL)\in G\T G^0/L;x\i gx\in\Si^*\}=\cup_{j\in J}\tY_{L,S_j}$$
(a disjoint union),
$$Y_{L,\Si}=\cup_{x\in G^0}x\Si^*x\i=\cup_{j\in J}Y_{L,S_j}.$$
(a not necessarily disjoint union). Thus $Y_{L,\Si}$ is a finite union of strata of 
equal dimension (see 3.13(a)) of $G$. Hence $Y_{L,\Si}$ is a locally closed smooth
subvariety of pure dimension of $G$ and any of its irreducible components is of the 
form $Y_{L,S_j}$ for some $j$ (which is not necessarily unique). Let 
$\bY_{L,S_j},\bY_{L,\Si}$ be the closure of $Y_{L,S_j},Y_{L,\Si}$ in $G$. We define
a local system $\tce$ on $\tY_{L,\Si}$ by the requirement that $b^*\ce=a^*\tce$ 
where $a(g,x)=(g,xL),b(g,x)=x\i gx$ in the diagram
$$\tY_{L,\Si}@<a<<\{(g,x)\in G\T G^0;x\i gx\in\Si^*\}@>b>>\Si.$$
(We use the fact that $a$ is a principal $L$-bundle and $b^*\ce$ is 
$L$-equivariant.)

Define $\p:\tY_{L,\Si}@>>>Y_{L,\Si}$ by $\p(g,xL)=g$. Using 3.13(a), we see that 
$\p$ is a finite unramified covering (for any irreducible component $U$ of 
$Y_{L,\Si}$, $\p:\p\i(U)@>>>U$ may be identified with 
$\sqc_{j\in J;Y_{L,S_j}=U}\tY_{L,S_j}@>>>U$ given by 
$\p_j:\tY_{L,S_j}@>>>Y_{L,S_j}$, $(g,xL)\m g$). It follows that $\p_!\tce$ is a 
local system on $Y_{L,\Si}$ such that for any irreducible component $U$ of 
$Y_{L,\Si}$, we have
$$\p_!\tce|_U=\op_{j\in J;Y_{L,S_j}=U}\p_{j!}\tce_j$$
where $\tce_j$ is the local system on $\tY_{L,S_j}$ defined in terms of $\ce_j$ as 
in 5.6. Then $IC(\bY_{L,\Si},\p_!\tce)\in\cd(\bY_{L,\Si})$ is well defined and 
$$IC(\bY_{L,\Si},\p_!\tce)=\op_{j\in J}IC(\bY_{L,S_j},\p_{j!}\tce_j)\tag a$$
where $IC(\bY_{L,S_j},\p_{j!}\tce_j)$ is regarded as a complex on $\bY_{L,\Si}$ 
which is zero on $\bY_{L,\Si}-\bY_{L,S_j}$.

\proclaim{Proposition 15.2}Assume that 
$S_j^\o=\{g\in S_j;g\text{ unipotent}\}\ne\em$ for any $j\in J$. Let 
$\Si^\o=\cup_{j\in J}S_j^\o$, let $p:\Si@>>>\Si$ be the morphism $p(g)=g_u$ and let
$\ce^1=p^*\ce$. Then $\ce^1|_{S_j}\in\cs(S_j)$ for all $j\in J$ and we can define 
$\tce^1$ in terms of $\ce^1$ just like $\tce$ is defined in terms of $\ce$. Let 
$\bY_{L,\Si}^\o=\{g\in\bY_{L,\Si};g\text{ unipotent}\}$. We have a canonical 
isomorphism

(a) $IC(\bY_{L,\Si},\p_!\tce)|_{\bY_{L,\Si}^\o}@>\si>>
IC(\bY_{L,\Si},\p_!\tce^1)|_{\bY_{L,\Si}^\o}$
\nl
in $\cd(\bY_{L,\Si}^\o)$.
\endproclaim
In other words, $IC(\bY_{L,\Si},\p_!\tce)|_{\bY_{L,\Si}^\o}$ depends only on
$\ce|_{\Si^\o}$, not on $\ce$ on the whole of $\Si$. The proof will be given in 
15.5-15.11.

\subhead 15.3\endsubhead
We return to the setup in 15.1. Let $P$ be a parabolic of $G^0$ with Levi $L$ such 
that $\Si\sub N_GP$. By 3.14 we have $\bY_{L,S_j}=\cup_{x\in G^0}x\bS_jU_Px\i$ hence
$$\bY_{L,\Si}=\cup_{j\in J}\bY_{L,S_j}=\cup_{x\in G^0}x\baSi U_Px\i.$$
Let $X_j=\{(g,xP)\in G\T G^0/P;x\i gx\in\bS_jU_P\}$. Let 

$X=\{(g,xP)\in G\T G^0/P;x\i gx\in\baSi U_P\}=\cup_{j\in J}X_j$.
\nl
Define $\ps:X@>>>\bY_{L,\Si}$ by $\ps(g,xP)=g$. This map is onto. We have the 
following generalization of Lemma 5.5:

(a) {\it $(g,xL)\m(g,xP)$ is an isomorphism $\tY_{L,\Si}@>\si>>\ps\i(Y_{L,\Si})$.}
\nl
We verify this only at the level of sets. The proof of injectivity is the same as 
that in Lemma 5.5. We prove surjectivity. Let $(g,xP)\in\ps\i(Y_{L,\Si})$. We have 
$(g,xP)\in X_j$ for some $j\in J$ (hence $g\in\bY_{L,S_j}$) and $g\in Y_{L,S_{j'}}$
for some $j'\in J$. Since $Y_{L,S_{j'}}$ is a stratum of $G$ that meets 
$\bY_{L,S_j}$ (which is a union of strata, see 3.15) we see that 
$Y_{L,S_{j'}}\sub\bY_{L,S_j}$. Since $\dim Y_{L,S_{j'}}=\dim Y_{L,S_j}$ (see 
3.13(b)) and $Y_{L,S_j}$ is the only stratum in $\bY_{L,S_j}$ of its dimension, we 
see that $Y_{L,S_{j'}}=Y_{L,S_j}$. Thus, $g\in Y_{L,S_j}$. Using the surjectivity of
the map in Lemma 5.5 (for $S_j$ instead of $S$) we see that there exists 
$(g,x'L)\in\tY_{L,S_j}$ such that $(g,x'P)=(g,xP)$. This proves the surjectivity of
our map.

\subhead 15.4\endsubhead
We are still in the setup of 15.1. For any stratum $S'$ of $N_LG$ that is contained 
in $\baSi$ let $X_{S'}=\{(g,xP)\in G\T G^0/P;x\i gx\in S'U_P\}\sub X$. We have 
$X=\sqc_{S'}X_{S'}$ (a finite union over all $S'$ as above). Let 
$X_\Si=\cup_{j\in J}X_{S_j}$. Then $X_\Si$ is a smooth, open dense subvariety of 
pure dimension of $X$ (since $\Si$ is a smooth, open dense subvariety of pure 
dimension of $\baSi$).

We define a local system $\bce$ on $X_\Si$ by the requirement that 
$b'{}^*\ce=a'{}^*\bce$ where $a'(g,x)=(g,xP),b'(g,x)=f(x\i gx)$ ($f$ as in 5.4) in 
the diagram
$$X_\Si@<a'<<\{(g,x)\in G\T G^0;x\i gx\in\Si U_P\}@>b'>>\Si.$$
(We use the fact that $a'$ is a principal $P$-bundle and $b'{}^*\ce$ is 
$P$-equivariant.) Then $IC(X,\bce)\in\cd(X)$ is well defined and we have
$a''{}^*IC(X,\bce)=b''{}^*IC(\baSi,\ce)$ where $a''(g,x)=(g,xP),b''(g,x)=f(x\i gx)$
in the diagram
$$X@<a''<<\{(g,x)\in G\T G^0;x\i gx\in\baSi U_P\}@>b''>>\baSi.$$
We have the following generalization of Lemma 5.7:

(a) {\it $\ps_!(IC(X,\bce))$ is canonically isomorphic to 
$IC(\bY_{L,\Si},\p_!\tce)$.}
\nl
The proof is similar to that of Lemma 5.7. Let $K=IC(X,\bce)$ and let 
$K^*=IC(X,\bce^*)$ where $\bce^*$ is defined like $\bce$, replacing $\ce$ by the 
dual local system. Using 15.3(a), we see that $\ps_!K|_{Y_{L,\Si}}=\p_!\tce$. As in 
the proof of Lemma 5.7, it is enough to verify the following statement.

For any $i>0$ we have $\dim\supp\ch^i(\ps_!K)<\dim\bY_{L,\Si}$ and
$\dim\supp\ch^i(\ps_!K^*)<\dim\bY_{L,\Si}-i$.
\nl
We shall verify this only for $K$; the corresponding statement for $K^*$ is entirely
similar. As in the proof of Lemma 5.7, it is enough to prove:

For any $i>0$ and any stratum $S'$ of $N_GL$ contained in $\baSi$ we have

$\dim\{g\in\bY_{L,\Si};H^i_c(\ps\i(g)\cap X_{S'},K)\ne 0\}<\dim\bY_{L,\Si}-i$.
\nl
Assume first that $S'$ is not one of the $S_j$. As in the proof of Lemma 5.7, we see
that it is enough to prove:
$$\dim\{g\in\bY_{L,\Si};\dim(\ps\i(g)\cap X_{S'})
>\fra{i}{2}-\fra{1}{2}(\dim\Si-\dim S')\}<\dim\bY_{L,\Si}-i.$$
Since $\bY_{L,\Si}=\cup_{j\in J}\bY_{L,S_j}$ this follows from the analogous
inequality in the proof of Lemma 5.7 where $\bY_{L,\Si}$ is replaced by 
$\bY_{L,S_j}$.

Next assume that $S'=S_j$ for some $j\in J$. As in the proof of Lemma 5.7, we see 
that it is enough to prove:
$$\dim\{g\in\bY_{L,\Si};\dim(\ps\i(g)\cap X_{S_j})\ge\fra{i}{2}\}<\dim\bY_{L,\Si}-i.
$$
If $g\in\bY_{L,\Si}$ satisfies $\dim(\ps\i(g)\cap X_{S_j})\ge\fra{i}{2}$ then 
$\ps\i(g)\cap X_{S_j}\ne\em$ and we have $x\i gx\in S_jU_P$ for some $x\in G^0$.
Hence $g\in\ps(X_j)=\bY_{L,S_j}$. Thus it is enough to prove:
$$\dim\{g\in\bY_{L,S_j};\dim(\ps\i(g)\cap X_{S_j})\ge\fra{i}{2}\}<\dim\bY_{L,S_j}-i.
$$
This is actually contained in the proof of Lemma 5.7. This completes the proof.

\subhead 15.5\endsubhead
The remainder of this section (except 15.12) is concerned with the proof of 
Proposition 15.2. If the analogue of the isomorphism 15.2(a) is known when $\Si$ is
replaced by $S_j$ then, taking direct sum over $j\in J$ and using the isomorphism 
15.1(a) and its analogue for $\ce^1$ instead of $\ce$, we obtain an isomorphism as 
in 15.2(a) for $\Si$. Thus, to prove Proposition 15.2, we may assume that $\Si=S$ is
a single isolated stratum of $N_GL$ with $(L,S)\in\AA$. Then $\ce\in\cs(S)$. 

Let $D$ be the connected component of $G$ that contains $S$. Let $\d$ be the 
connected component of $N_GL$ that contains $S$. By the assumption of Proposition 
15.2, the set $S^\o$ of unipotent elements in $S$ is a single $L$-conjugacy class. 
We have $S={}^\d\cz^0_LS^\o$. Let 
$S^*,Y=Y_{L,S},\bY=\bY_{L,S},\tY=\tY_{L,S},\p:\tY@>>>Y$ be as in 15.1. 

Let $\mu_\iy=\mu_\iy({}^\d\cz_L^0)$ be as in 5.3. There is a canonical direct sum 
decomposition $\ce=\op_\l\ce^\l$ in $\cs(S)$ where $\l$ runs over the set of 
homomorphisms $\mu_\iy@>>>\bbq^*$ that factor through some 
$\mu_n=\mu_n({}^\d\cz_L^0)$ (see 5.3). Here $\ce^\l$ has the property that, for any
$z\in{}^\d\cz_L^0,h\in S^\o$, the monodromy action of $\mu_\iy$ on the stalk 
$\ce^l_{zh}$ of the local system $\ce^\l|_{{}^\d\cz_L^0h}$ (equivariant for the 
transitive ${}^\d\cz_L^0$-action $z_1:zh\m z_1^nzh$ for some $n$ invertible in 
$\kk$) is through $\l$. We have canonically 
$$IC(\bY,\p_!\tce)=\op_\l IC(\bY,\p_!\ti{\ce^\l}),\qua
                    IC(\bY,\p_!\tce^1)=\op_\l IC(\bY,\p_!\ti{\ce^{1\l}}).$$
Hence if we can construct the isomorphism 15.2(a) for each $\ce^\l$ instead of $\ce$
then by taking direct sums we get the isomorphism (a) for $\ce$. Thus we may assume
that $\ce=\ce^\l$ for some $\l$. 

Let $\cl$ be a local system of rank $1$ on ${}^\d\cz_L^0$, equivariant for the 
transitive ${}^\d\cz_L^0$-action $z_1:z\m z_1^nz$ for some $n$ invertible in $\kk$,
whose associate homomorphism $\mu_\iy@>>>\bbq^*$ is $\l\i$; let $v_0$ be a basis
vector of the stalk $\cl_1$. The pair $(\cl,v_0)$ is defined up to a unique 
isomorphism. 

Let $b:S@>>>{}^\d\cz_L^0$ be $g\m g_s$. From the definitions, the restriction of 
$\ce\ot b^*\cl$ to any fibre of $S@>>>S^\o,g\m g_u$ is isomorphic to $\bbq^e$ where 
$e$ is the rank of $\ce$. Also 
$$(\ce\ot b^*\cl)|_{S^\o}=\ce^1|_{S^\o}=\ce|_{S^\o}\tag a$$
canonically, using the identification $\cl_1=\bbq$, $v_0\lra 1$. It follows that 
$\ce\ot b^*\cl$ is isomorphic to $\ce^1$ and there is a unique isomorphisms between
these two local systems which induces for the restrictions to $S^\o$ the 
identification (a). Thus we have a canonical isomorphism

(b) $\ce\ot b^*\cl\cong\ce^1$.
\nl
It is then enough to construct a canonical isomorphism
$$IC(\bY,\p_!\tce)|_{\bY^\o}@>\si>>IC(\bY,\p_!\wt{\ce\ot b^*\cl})|_{\bY^\o}.\tag c$$
where $\bY^\o=\{g\in\bY;g\text{ unipotent}\}$. 

\subhead 15.6\endsubhead
Let $S^\da$ be the subset of $S$ consisting of those 
$zh,z\in{}^\d\cz^0_L,h\in S^\o$ such that 
$$n\in N_{G^0}\d=\{x\in G^0;x\d x\i=\d\},nz=zn\imp n\in L$$
(that is, $z$ has trivial stabilizer for the conjugation action of $N_{G^0}\d/L$ on 
${}^\d\cz^0_L$). Now $S^\da$ is open dense in $S$. (Since $N_{G^0}\d/L$ is finite, 
it is enough to show that, for any $n\in N_{G^0}\d-L$, the closed subset 
$\{zh;z\in{}^\d\cz^0_L,h\in S^\o,nz=zn\}$ of $S$ is $\ne S$ or that the closed 
subset $\{z\in{}^\d\cz^0_L;nz=zn\}$ of ${}^\d\cz^0_L$ is $\ne{}^\d\cz^0_L$. If it is
equal to ${}^\d\cz^0_L$, then $n\in Z_{G^0}({}^\d\cz^0_L)$ hence using 1.10(a), 
$n\in L$, a contradiction.) 

Since $S^*$ is an open dense subset of $S$ it follows that $S^{*\da}=S^*\cap S^\da$
is an open dense subset of $S$. Hence
$$\tY^\da=\{(g,xL)\in G\T G^0;x\i gx\in S^{*\da}\}$$
is an open dense subset of $\tY$ (we use that $S^{*\da}$ is stable under 
$L$-conjugacy). Also, $\tY^\da$ is stable under the free action of $\cw_S$ on $\tY$,
see 3.13; it follows that $Y^\da=\p(\tY^\da)$ is open dense in $Y$ and 
$\tY^\da=\p\i(Y^\da)$.

Let $\s:D@>>>D//G^0$ be as in 7.1. Let $A=\s(\bY)$. Let $u$ be a unipotent, 
quasi-semisimple element of $N_GL$ such that $u\in\d$. As in the proof of 7.3(b) we
see that $A=\{\s(zu);z\in{}^\d\cz_L^0\}$. Let 
$$\Up=\{(g,z)\in\bY\T{}^\d\cz_L^0;\s(g)=\s(zu)\}.$$
Define $\k:\tY^\da@>>>\Up$ by $(g,xL)\m(g,z)$ where 
$x\i gx=zh\in S^{*\da},z\in{}^\d\cz^0_L,h\in S^\o$. (This definition is correct 
since $\s(zh)=\s(zu)$ for $z\in{}^\d\cz^0_L,h\in S^\o$, see 7.3.) Let $\Up'$ be the
closure of $\k(\tY^\da)$ in $\Up$. Define $\ps':\Up'@>>>\bY$ by $(g,z)\m g$. 

\proclaim{Lemma 15.7}(a) $\k(\tY^\da)$ is open in $\Up'$.  

(b) $\k$ restricts to an isomorphism $\tY^\da@>\si>>\k(\tY^\da)$.

(c) $\ps'$ is a finite surjective morphism.
\endproclaim
We verify (b) at the level of sets. We must show that $\k:\tY^\da@>>>\Up$ is 
injective. Assume that $(g,xL),(g',x'L)$ in $\tY^\da$ have the same image under 
$\k$. Then $g=g'$ and 
$$x\i gx=zh,x'{}\i gx'=zh',z\in{}^\d\cz^0_L,h\in S^\o,h'\in S^\o,zh\in S^{*\da}.$$
From 3.13(a) we see that $x'=xn\i$ for some $n\in N_{G^0}L,nS^\o n\i=S^\o$. Then 
$zh'=nzhn\i=nzn\i(nhn\i)$. Now $z,z'$ are semisimple elements commuting with the 
unipotent elements $h',nhn\i$. It follows that $z=nzn\i$. Since $n\in N_{G^0}\d$ 
and $z$ has trivial stabilizer in $N_{G^0}\d/L$, we see that $n\in L$. Thus 
$xL=x'L$ and $\k$ is injective, as required.

We prove (a). Let $P$ be a parabolic of $G^0$ with Levi $L$ such that $S\sub N_GP$.
Let $\ps:X@>>>\bY$ be as in 3.14. Define $f:X@>>>\Up$ by $f(g,xP)=(g,z)$ where 
$x\i gx\in z\bS^\o U_P,z\in{}^\d\cz^0_L$. (We show that $z$ is uniquely determined 
by $(g,xP)$. It is enough to show that, if $zhv=p\i z'h'v'p$ with 
$z,z'\in{}^\d\cz_L^0,h,h'\in\bS^\o,v,v'\in U_P,p\in P$, then $z=z'$. Writing 
$p\in lU_P$ with $l\in L$, we have $zh=l\i z'h'l=z'l\i h'l$. By the uniqueness of 
Jordan decomposition we have $z=z'$ as required. To show that $f$ is well defined we
must also show that, if $g\in\bY,x\in G^0,z\in{}^\d\cz_L^0,h\in\bS^\o,v\in U_P$ 
satisfy $x\i gx=zhv$ then $\s(zu)=\s(g)$. Clearly, $\s(g)=\s(zhv)$ hence we must 
show that $\s(zu)=\s(zhv)$. This follows from the description of $\s$ given in 7.1.)
We have $\ps=pr_1\circ f$ where $pr_1:\Up@>>>\bY$ is $(g,z)\m z$. Since $\ps$ is 
proper, it follows that $f$ is proper and $f(X)$ is closed in $\Up$. Since $Y^\da$ 
is open in $Y$ and $Y$ is open in $\bY$ we see that $Y^\da$ is open in $\bY$ and 
$\ps\i(\bY-Y^\da)$ is closed in $X$. Since $f$ is proper it follows that 
$f(\ps\i(\bY-Y^\da))$ is closed in $f(X)$. We have
$f(\ps\i(\bY-Y^\da))\cap f(\ps\i(Y^\da))=\em$. (Indeed if $(g,xP),(g',x'P)\in X$ 
have the same image under $f$ and $g\in Y^\da$, then $g'\in Y^\da$.) Thus, 
$f(\ps\i(Y^\da))$ is the complement of $f(\ps\i(\bY-Y^\da))$ in $f(X)$ hence it is 
open in $f(X)$. From Lemma 5.5 we see that 
$\g_0:\tY^\da@>>>\ps\i(Y^\da),(g,xL)\m(g,xP)$ is an isomorphism. Since $f\g_0=\k$, 
we see that $\k(\tY^\da)=f(\ps\i(Y^\da))$ hence $\k(\tY^\da)$ is open in $f(X)$. 
Since $\Up'$ is the closure of $\k(\tY^\da)$ in $\Up$ and $f(X)$ is a closed subset
of $\Up$ containing $\k(\tY^\da)$, we see that $\Up'\sub f(X)$. Since $X,\bY$ are 
irreducible of the same dimension (see the proof of 3.14) and $\tY^\da$ is open 
dense in $\bY$, we see that $f(X),\k(\tY^\da),\Up'$ are irreducible and 
$$\dim\bY=\dim\tY^\da=\dim\k(\tY^\da)\le\dim\Up'\le\dim f(X)\le\dim X=\dim\bY.$$
It follows that $\dim\Up'=\dim f(X)=\dim X$ hence $\Up'=f(X)$. Thus, $\k(\tY^\da)$
is open in $\Up'$.

We prove (c). Let $T_1$ be a maximal torus of $Z_G(u)^0$ that contains 
${}^\d\cz^0_L$. As in 7.1, $\s$ induces a finite morphism $uT_1@>>>D//G^0$. This 
restricts to a finite morphism $u{}^\d\cz^0_L@>>>A$. Since $pr_1:\Up@>>>\bY$ is 
obtained from this finite morphism by change of base, it follows that $pr_1$ is a 
finite morphism. Restricting to the closed subset $\Up'$ of $\Up$, we deduce that 
$\ps':\Up'@>>>\bY$ is a finite morphism. To see that it is surjective, we note that 
$\ps:X@>>>\bY$ is surjective and $\ps$ factorizes as $X@>f>>\Up'@>\ps'>>\bY$. The 
lemma is proved.

\subhead 15.8\endsubhead
Let $\uce=\k_!(\tce|_{\tY^\da})$, a local system on $\k(\tY^\da)$. Since $\tY$ is
smooth (see 3.17) and $\tY^\da$ is open in $\tY$, we see that $\tY^\da$ is smooth. 
Using 15.7(a),(b) we see that $\k(\tY^\da)$ is a smooth open dense subvariety of 
$\Up'$. Hence $IC(\Up',\uce)\in\cd(\Up')$ is well defined.

\proclaim{Lemma 15.9} $\ps'_!IC(\Up',\uce)$ is canonically isomorphic to 
$IC(\bY,\p_!\tce)$.
\endproclaim
Since $Y^\da$ is open dense in $Y$, we have canonically
$IC(\bY,\p_!\tce)=IC(\bY,(\p_!\tce)_{Y^\da})$.

Let $K_1=IC(\Up',\uce),K^*_1=IC(\Up',\uce^*)$ where $\uce^*$ is defined like $\uce$
replacing $\ce$ by $\ce^*$. Then $K^*_1$ is the Verdier dual of $K_1$ with a 
suitable shift. Since $\ps'$ is proper, it follows that $\ps'_!(K^*_1)$ is the 
Verdier dual of $\ps'_!K_1$ with a suitable shift. From the definitions it is clear
that $(\ps'_!K_1)|_{Y^\da}=(\p_!\tce)|_{Y^\da}$. By the definition of an 
intersection cohomology complex we see that it is enough to verify the following 
statement.

For any $i>0$ we have 

$\dim\supp\ch^i(\ps'_!K_1)<\dim\bY-i$ and $\dim\supp\ch^i(\ps'_!(K^*_1))<\dim\bY-i$.
\nl
We shall only verify this for $K_1$; the corresponding statement for $K^*_1$ is 
entirely analogous. For $g\in\bY$, $\ps'{}\i(g)$ is finite hence 
$$\ch^i_g(\ps'_!K_1)=H^i_c(\ps'{}\i(g),K_1)=\op_{y\in\ps'{}\i(g)}\ch^i_yK_1.$$
It is enough to check that $\dim\ps'(\{y\in\Up';\ch^i_yK_1\ne 0\})<\dim\bY-i$. But 
$$\dim\ps'(\{y\in\Up';\ch^i_yK_1\ne 0\})\le\dim\{y\in\Up';\ch^i_yK_1\ne 0\}<
\dim\Up'-i=\dim\bY-i$$
(since $K_1$ is an intersection cohomology complex on $\Up'$). The lemma is proved.

\proclaim{Lemma 15.10} Let $\cl'=pr_2^*\cl$ where $pr_2:\Up'@>>>{}^\d\cz_L^0$ is the
second projection. The restriction of the local system $\cl'$ on $\Up'$ to 
$\k(\tY^\da)$ is denoted again by $\cl'$. We have a canonical isomorphism

(a) $IC(\Up',\uce)\ot\cl'@>\si>>IC(\Up',\uce\ot\cl')$.
\endproclaim
The restrictions of the two sides of (a) to the open dense subset $\k(\tY^\da)$ of 
$\Up'$ are canonically isomorphic (they can both be identified with $\uce\ot\cl'$).
From the properties of intersection cohomology complexes it follows that this
extends to an isomorphism as in (a) provided we can show that the left hand side of
(a) is an intersection cohomology complex on $\Up'$. To do this we choose a
parabolic $P$ as in the proof of 15.7(a). Let $f:X@>>>\Up'$ be as in that proof. 
Then $f$ is proper, surjective. Let $X_S,\bce$ be as in 5.6 and let $K,K^*$ be as in
5.7. We show that

(b) $f_!K=IC(\bY',\uce)$.
\nl
It is clear that $(f_!K)|_{\k(\tY^\da)}=\uce$. As in the proof of Lemma 5.7, it is 
enough to verify the following statement.

For any $i>0$ we have 

$\dim\supp\ch^i(f_!K)<\dim\Up'-i$ and $\dim\supp\ch^i(f_!K')<\dim\Up'-i$.
\nl
We shall only verify this for $K$; the corresponding statement for $K^*$ is entirely
similar. Let $y\in\Up'$. If $\ch^i_y(f_!K)\ne 0$ then
$\op_{y'\in\Up';\ps'(y')=\ps'(y)}\ch^i_{y'}(f_!K)\ne 0$ hence 
$\ch^i_{\ps'(y)}(\ps'_!f_!K)\ne 0$. (We use that $\ps'{}\i(\ps'(y))$ is finite.) 
Thus, $y\in\supp\ch^i(f_!K)$ implies $\ps'(y)\in\supp\ch^i(\ps_!K)$ that
$y\in\ps'{}\i(\supp\ch^i(\ps_!K))$. We see that 

$\supp\ch^i(f_!K)\sub\ps'{}\i(\supp\ch^i(\ps_!K))$. 
\nl
Hence 

$\dim\supp\ch^i(f_!K)\le\dim\ps'{}\i(\supp\ch^i(\ps_!K))$.
\nl
Since $\ps'$ is finite, surjective, we have 
$\dim\ps'{}\i(\supp\ch^i(\ps_!K))\le\dim\supp\ch^i(\ps_!K)$. It is then enough to 
show that $\dim\supp\ch^i(f_!K)<\dim\Up'-i$. But this follows from an estimate in 
the proof of Lemma 5.7 since $\dim\Up'=\dim Y$. This proves (b).

We show that 

(c) $K\ot f^*\cl'=IC(X,\bce^1)$
\nl
where $\bce^1$ is the local system on $X_S$ defined in terms of $\ce^1$ in the same
way as $\bce$ is defined in terms of $\ce$ (see 5.6). With the notation in 5.6(a), 
it is enough to show that $a''{}^*(K\ot f^*\cl')=a''{}^*IC(X,\bce^1)$ or that
$b''{}^*IC(\bS,\ce)\ot b''{}^*\bar b^*\cl=b''{}^*IC(\bS,\ce^1)$ where 
$\bar b:\bS@>>>{}^\d\cz_L^0$ is $g\m g_s$. It is enough to show that 
$IC(\bS,\ce)\ot\bar b^*\cl=IC(\bS,\ce^1)$. This follows immediately from the
definitions. Thus (c) is proved.

Using (b),(c), we have
$$IC(\Up',\uce)\ot\cl'=f_!K\ot\cl'=f_!(K\ot f^*\cl')=f_!(IC(X,\bce^1))
=IC(\Up',\uce^1),$$
the last step being (b) applied to $\ce^1$ instead of $\ce$. (We define $\uce^1$ in
terms of $\ce^1$ in the same way as $\uce$ is defined in terms of $\ce$.) This 
completes the proof.

It is clear that, although $P$ is used in the proof above, the isomorphism (a) that
we construct does not depend on the choice of $P$ hence it is truly canonical.

\subhead 15.11\endsubhead
Let $\Upo=\bY^\o\T\{1\}\sub\bY\T{}^\d\cz_L^0$. Then $\Upo$ is a closed subset of
$\Up$. More precisely, we have $\Upo\sub\Up'$. (With notation in the proof of Lemma
15.7, it is enough to show that $\Upo\sub f(X)$. Let $g\in\bY^\o$. Since 
$\ps:X@>>>\bY$ is surjective, we can find $xP\in G^0/P$ such that 
$x\i gx\in z\bS^\o U_P$ where $z\in{}^\d\cz_L^0$. Since $x\i gx$ is unipotent, we 
must have $z=1$. Hence $f(g,xP)=(g,1)$. Thus $(g,1)\in f(X)$ as required.) Thus we 
have $\Upo\sub\ps'{}\i(\bY^\o)$. This in fact an equality. (If $(g,z)\in\Up'$ and 
$g$ is unipotent, then from $\s(zu)=\s(g)=\o$ we see that $zu$ is unipotent hence 
$z=1$.) Since $\ps'$ restricts to an isomorphism $\Upo@>\si>>\bY^\o$, it also 
restricts to an isomorphism $\ps'_0:\ps'{}\i(\bY^\o)@>\si>>\bY^\o$. Via $\ps'_0$ we
may identify $\bY^\o$ with $\Upo=\ps'{}\i(\bY^\o)$. By change of base, we have
$$IC(\bY,\p_!\tce)|_{\bY^\o}=(\ps'_!IC(\Up',\uce))|_{\bY^\o}=
\ps'_{0!}(IC(\Up',\uce)|_{\ps'{}\i(\bY^\o)}).$$
Thus, via $\ps'_0$, we may identify 
$IC(\bY,\p_!\tce)|_{\bY^\o}=IC(\Up',\uce)|_{\Upo}$.

Similarly, we may identify 
$IC(\bY,\p_!\wt{\ce\ot b^*\cl})|_{\bY^\o}=IC(\Up',\un{\ce\ot b^*\cl})|_{\Upo}$. 
Hence, in order to construct the isomorphism 15.5(c), it is enough to construct a 
canonical isomorphism

$IC(\Up',\uce)|_{\Upo}@>\si>>IC(\Up',\un{\ce\ot b^*\cl})|_{\Upo}$.
\nl
From the definitions we have $\un{\ce\ot b^*\cl}=\uce\ot\cl'$. Hence it is enough to
construct a canonical isomorphism

$IC(\Up',\uce)|_{\Upo}@>\si>>IC(\Up',\uce\ot\cl')|_{\Upo}$.
\nl
We have canonically $\cl'_{\Upo}=\bbq$ (we identify $\cl_1=\bbq$ by $v_0\lra 1$).
Hence $(IC(\Up',\uce)\ot\cl')|_{\Upo}=IC(\Up',\uce)|_{\Upo}$. Thus it is enough to
construct a canonical isomorphism

$(IC(\Up',\uce)\ot\cl')|_{\Upo}@>\si>>IC(\Up',\uce\ot\cl')|_{\Upo}$.
\nl
This is obtained by restricting to $\Upo$ the isomorphism 15.10(a). This completes 
the proof of Proposition 15.2.

\subhead 15.12\endsubhead
In the remainder of this section we assume that $\kk$ is an algebraic closure of a 
finite field $\FF_q$ and that $G$ has a fixed $\FF_q$-structure with Frobenius map 
$F:G@>>>G$. 

For any algebraic variety $Z$ defined over $\FF_q$ with Frobenius map $F:Z@>>>Z$, an
object $A\in\cd(Z)$ and an isomorphism $\ph:F^*A@>\si>>A$ in $\cd(Z)$ we define the
{\it characteristic function} $\c_{A,\ph}:Z^F@>>>\bbq$ by
$$\c_{A,\ph}(z)=\sum_i(-1)^i\tr(\ph,\ch^i_zA),\qua (z\in Z^F).\tag a$$
(The map induced by $\ph$ on the stalk $\ch^i_zA$ is denoted again by $\ph$.)

Consider a quadruple $(L,\Si^\o,\cf,\ph_1)$ where 

$L$ is an $F$-stable Levi of some (not necessarily $F$-stable) parabolic of $G^0$;

$\Si^\o$ is the set of unipotent elements in a subset $\Si$ of $N_GL$ as in 15.1,
15.2 such that $F(\Si^\o)=\Si^\o$ (or equivalently $F(\Si)=\Si$); note that 
$\Si^\o,\Si$ determine each other;

$\cf$ is an $L$-equivariant local system on $\Si^\o$;

$\ph_1:F^*\cf@>\si>>\cf$ is an isomorphism of local systems on $\Si^\o$. 
\nl
Consider the complex $\fK=IC(\bY_{L,\Si},\p_!\tce)\in\cd(\bY_{L,\Si})$ where 
$\bY_{L,\Si},\p:\tY_{L,\Si}@>>>Y_{L,\Si}$ are defined as in 15.1 and $\tce$ is 
defined as in 15.1 in terms of a local system $\ce$ on $\Si$ (as in 15.1) such that
$\cf$ is the inverse image of $\ce$ under the inclusion $\Si^\o@>>>\Si$. We assume 
that we are given an isomorphism $\ph'_1:F^*\ce@>\si>>\ce$ of local systems on $\Si$
extending $\ph_1$. (Note that we can always find $\ce,\ph'_1$ as above: for example,
we have the "trivial choice" where $\ce$ is the inverse image of $\cf$ under 
$\Si@>>>\Si^\o,g\m g_u$ and $\ph'_1$ is induced by $\ph_1$. However, for future 
applications, it is necessary to allow other choices of $\ce,\ph'_1$.) Now 
$\bY_{L,\Si},\tY_{L,\Si},Y_{L,\Si}$ have natural $\FF_q$-structures with Frobenius 
maps $F$ and $\ph'_1$ induces an isomorphism $F^*\tce@>\si>>\tce$ of local systems 
on $\tY_{L,\Si}$, an isomorphism $F^*\p_!\tce@>\si>>\p_!\tce$ of local systems on 
$Y_{L,\Si}$ and an isomorphism $\ph:F^*\fK@>\si>>\fK$ in $\cd(\bY_{L,\Si})$. We 
define a function
$$Q_{L,G,\Si^\o,\cf,\ph_1}:\{\text{unipotent elements in } G^F\}@>>>\bbq\tag b$$
by
$$Q_{L,G,\Si^\o,\cf,\ph_1}(u)=\c_{\fK,\ph}(u)\tag c$$
(see (a)) if $u\in \bY_{L,\Si}$ and $Q_{L,G,\Si^\o,\cf,\ph_1}(u)=0$ if 
$u\notin\bY_{L,\Si}$. The function (c) is called a {\it generalized Green function}.
It extends (up to a sign) a definition given in \cite{\CSII, 8.3.1} (in the case 
where $G=G^0$ and $\Si^\o$ is a single unipotent class). From Proposition 15.2 we 
see that 

(d) {\it $Q_{L,G,\Si^\o,\cf,\ph_1}(u)$ is independent of the choice of 
$\ce,\ph'_1$},
\nl
namely it is the same for a general $\ce,\ph'_1$ as for the "trivial choice". (The 
isomorphism in 15.2 is compatible with the Frobenius maps.)

\head 16. The characteristic function $\c_{\fK,\ph}$\endhead
\subhead 16.1\endsubhead
In this section we fix $(L,S)\in\AA$ and $\ce\in\cs(S)$. Let $\d$ be the connected 
component of $N_GL$ that contains $S$. Recall (cf. 1.22) that 

(a) $S_s=\{g_s;g\in S\}$ is a single ${}^\d\cz_L^0\T L$-orbit on $N_GL$ for the 
action $(z,x):y\m xzyx\i$.

\proclaim{Lemma 16.2} (a) Let $s'\in S_s$. Let 
$\boc_{s'}=\{v\in Z_G(s');v\text{ unipotent, }s'v\in S\}$. Then 
$\boc_{s'}=\sqc_{j\in J}\fc_j$ where $J$ is finite and $\fc_j$ are (unipotent)
$Z_L(s')^0$-conjugacy classes of dimension independent of $j$. 

(b) For any $j\in J$, the stratum $S_j$ of $Z_{N_GL}(s')$ that contains $\fc_j$ is 
${}^\d\cz_L^0\fc_j$. In particular, $\dim S_j$ is independent of $j$.
\endproclaim
First we note that 

(c) {\it any $G^0$-conjugacy class in $G$ has finite intersection with 
${}^\d\cz^0_Ls'$.}
\nl
This follows from 1.14(a),(d) applied to $g=s'$ and to a maximal torus of 
$Z_G(g)^0$ that contains ${}^\d\cz^0_L$. 

In particular, any $L$-conjugacy class in $N_GL$ has finite intersection with 
${}^\d\cz^0_Ls'$. Hence the group $\tZ=\{y\in L;ys'y\i\in{}^\d\cz^0_Ls'\}$ contains
$Z_L(s')$ as a subgroup of finite index. Thus, $\tZ^0=Z_L(s')^0$. From 16.1(a) we 
see that $\tZ$ acts transitively (by conjugation) on $\boc_{s'}$. Since $Z_L(s')^0$
is normal in $\tZ$, it follows that $\tZ$ permutes transitively the 
$Z_L(s')^0$-orbits in $\boc_{s'}$; hence all these orbits have the same dimension. 
This proves (a).

We prove (b). We have 
$S_j=\{z\in\cz_{Z_L(s')^0}^0;zv=vz\}^0\fc_j=T_{N_GL}(s'v)\fc_j$ where $v$ is any 
element of $\fc_j$. Since $s'v\in S$, $s'v$ is isolated in $N_GL$ hence 
$T_{N_GL}(s'v)={}^\d\cz_L^0$, see 2.2. This proves (b). 

\subhead 16.3\endsubhead
We fix a semisimple element $s\in G$ and a unipotent element $u\in Z_G(s)$ such that
$su\in\bY$, the closure of $Y_{L,S}$ in $G$. Let $P$ be a parabolic of $G^0$ with 
Levi $L$ such that $S\sub N_GP$. Let 
$$M=\{x\in G^0;x\i sx\in S_s\},\hM=\{x\in G^0;x\i sx\in S_sU_P\}.$$
Let $\G$ be the set of orbits for the $Z_G(s)^0\T L$ action $(h,l):x\m hxl\i$ on 
$M$. We show:

(a) {\it $\G$ is finite.}
\nl
We may assume that $M\ne\em$. Let $x_0\in M$. Using 16.1(a), we see that it is 
enough to show that $\{x\in G^0;x\i sx\in{}^\d\cz^0_L(x_0\i sx_0)\}$ is a union of 
finitely many orbits under left translation by $Z_G(s)^0$, or equivalently, a union
of finitely many orbits under left translation by $Z_G(s)$, which contains 
$Z_G(s)^0$ with finite index. It is enough to note that any $G^0$-conjugacy class in
$G$ has finite intersection with ${}^\d\cz^0_L(x_0\i sx_0)$; see 16.2(c).

The group $\{n\in N_{G^0}L;nSn\i=S\}$ acts on $\G$ by $n:\co\m\co n\i$; this induces
an action of $\cw_S$ (see 3.13) on $\G$.

Let $\haG$ be the set of orbits for the $Z_G(s)^0\T P$ action $(h,p):x\m hxp\i$ on 
$\hM$. Any orbit $\co$ in $\G$ is contained in a unique orbit $\hco$ in $\haG$.

(b) {\it The map $\G@>>>\haG$, $\co\m\hco$, is a bijection.}
\nl
We show that our map is injective. Let $x,x'\in M$ be such that $x'=hxp\i$ for some 
$(h,p)\in Z_G(s)^0\T P$. We must show that $x,x'$ are in the same $Z_G(s)^0\T L$ 
orbit. Replacing $x$ by an element in the same $Z_G(s)^0\T L$ orbit we may assume 
that $x'=xu\i$ for some $u\in U_P$. Let $s''=x\i sx,s'=ux\i sxu\i$. Then $s',s''$ 
belong to $S_s$ hence to $N_LG\cap N_GP$ and 
$$s''s'{}\i=u\i(s'us'{}\i)\in (N_LG\cap N_GP)\cap U_P=\{1\}.$$
Thus, $u\i s'us'{}\i=1$ that is $u\in U_P\cap Z_G(s')=U_P\cap Z_G(s')^0$. (We use 
1.11.) Then $\z:=xux\i=xu\i u ux\i\in Z_G(xu\i s'ux\i)^0=Z_G(s)^0$ and 
$x'=xu\i=\z\i x$. Since $\z\in Z_G(s)^0$ we see that $x,x'$ are in the same
$Z_G(s)^0\T L$ orbit, as required.

We show that our map is surjective. Let $x\in\hM$. It is enough to show that for 
some $v\in U_P$ we have $xv\in M$. Now $x\i sx\in N_GP$ is semisimple. Hence, using
1.4(a), $x\i sx$ normalizes $vLv\i$ for some $v\in U_P$. Replacing $x$ by $xv$ we 
may assume that $x\i sx\in N_GP\cap N_GL$. We have $x\i sx=g'g''$ where 
$g'\in S_s,g''\in U_P$. Since $S_s\sub N_GP\cap N_GL$, we have 
$g'{}\i x\i sx\in(N_GP\cap N_GL)\cap U_P=\{1\}$, see 1.26. Thus, $x\i sx=g'\in S_s$
and $x\in M$. This completes the proof of (b).

Since $\G$ is finite, it follows that $\haG$ is finite.

The orbits of $Z_G(s)^0$ acting by left translation on 
$\{xP\in G^0/P;x\i sx\in S_sU_P\}$ are complete varieties. (Indeed, such an orbit is
of the form $Z_G(s)^0/(Z_G(s)^0\cap xPx\i)$ where $x\in G^0,x\i sx\in S_sU_P$ and it
is enough to show that $Z_G(s)^0\cap xPx\i$ is a parabolic of $Z_G(s)^0$, or
equivalently that $Z_G(x\i sx)^0\cap P$ is a parabolic of $Z_G(x\i sx)^0$. This 
follows from 1.12(a) since $x\i sx\in N_GP$.) Hence these orbits are closed. Since 
there are only finitely many such orbits (their number is $|\haG|=|\G|$), these 
orbits are also open.

\mpb

Let $\d_1$ be the connected component of $Z_G(s)$ that contains $u$. 

\proclaim{Lemma 16.4}There exists an open subset $\cu$ of $\d_1$ such that

(i) $\cu$ contains any unipotent element in $\d_1$;

(ii) $g\cu g\i=\cu$ for all $g\in Z_G(s)^0$;

(iii) for any $P$ as in 16.3 we have $h\in\cu,x\in G^0,x\i shx\in\bS U_P\imp 
x\i h_sx\in{}^\d\cz^0_LU_P,x\i sx\in S_sU_P$.
\endproclaim
Let $D$ be a connected component of $G$. A subset of $D$ is said to be {\it stable}
if it is a union of fibres of the map $\s:D@>>>D//G^0$ in 7.1. Let $g\in D$ be 
quasi-semisimple and let $T_1$ be a maximal torus of $Z_G(g)^0$. From 7.1(a) we 
deduce:

(a) {\it a stable subset $\car$ of $D$ is closed in $D$ if and only if 
$\car\cap gT_1$ is closed in $gT_1$.}
\nl
Next we show:

(b) {\it Assume that $D$ contains some unipotent elements. Let $\car_0$ be a subset
of $G^0$ which is a union of $G^0$-conjugacy classes such that the intersection of 
$\car_0$ with some/any maximal torus in $G^0$ is closed in that torus. Then 
$\car:=\{g\in D;g_s\in\car_0\}$ is a closed and stable subset of $D$.}
\nl
We show that $\car$ is stable. Let $y\in\car$ and let $y'\in D$ be such that 
$\s(y)=\s(y')$. We must show that $y'\in\car$. Let $v\in y_uZ_G(y_s)^0$ be 
unipotent, quasi-semisimple in $Z_G(y_s)$ and let $v'\in y'_uZ_G(y'_s)^0$ be 
unipotent, quasi-semisimple in $Z_G(y'_s)$. Then $y'_sv'=zy_svz\i$ for some 
$z\in G^0$ (see 7.1). It follows that $y'_s=zy_sz\i$. Since $y_s\in\car_0$ we see 
that $y'_s\in\car_0$ hence $y'\in\car$. We show that $\car$ is closed in $D$. Let 
$g,T_1$ be as above. Let $T'_1$ be a maximal torus of $G^0$ that contains $T_1$. 
Since $g$ is unipotent, we have 

$\car\cap gT_1=\{gt;t\in T_1\cap\car_0\}=g(T'_1\cap\car_0)\cap gT_1$.
\nl
This is closed in $gT_1$ since $T'_1\cap\car_0$ is closed in $G$. This proves (b).

(c) {\it Let $s'\in G$ be semisimple and let $g'\in Z_G(s')$ be such that 
$g'_s\in Z_G(s')^0$ and $Z_G(s'g'_s)\sub Z_G(s')$. If $s'g'$ is isolated in $G$, 
then $s'g'_u=g'_us'$ is isolated in $G$.} 
\nl
Let $T_1$ be a maximal torus of $Z_G(s')^0$ such that $g'_s\in T_1$. Then $T_1$ is 
also a maximal torus of $Z_G(s'g'_s)^0$. With the notation of 1.5 we have

$\Lie Z_G(s'g'_s)=\ft\op\op_{\a\in R'}\fg_\a,\qua
\Lie Z_G(s')=\ft\op\op_{\a\in R''}\fg_\a$
\nl
where $R'\sub R''$ are subsets of $R$ and $\ft=\Lie T_1$. The centre of 
$Z_G(s'g'_s)^0$ is $\{t\in T_1;\a(t)=1\frl\a\in R'\}$ and the centre of $Z_G(s')^0$
is $\{t\in T_1;\a(t)=1\frl\a\in R''\}$. Since $R'\sub R''$, the centre of 
$Z_G(s'g'_s)^0$ contains the centre of $Z_G(s')^0$. Hence 

$\cz_{Z_G(s')^0}^0\sub\cz_{Z_G(s'g'_s)^0}^0$ and
$\cz_{Z_G(s')^0}^0\cap Z_G(g_u)\sub\cz_{Z_G(s'g'_s)^0}^0\cap Z_G(g'_u)$.
\nl
It follows that $T(s'g'_u)\sub T(s'g'_sg'_u)$ (see 2.1). Since $s'g'_sg'_u$ is 
isolated in $G$ we have $T(s'g'_sg'_u)\sub\cz_{G^0}$. It follows that 
$T(s'g'_u)\sub\cz_{G^0}$, hence $s'g'_u$ is isolated in $G$. This proves (c).

(d) {\it Let $s'\in G$ be semisimple and let $g\in Z_G(s')$ be such that 
$g_s\in Z_G(s')^0$, $s'g\in\bS U_P$ ($P$ as in 16.3), $s'g_s\in S_s$ and
$Z_G(s'g_s)\sub Z_G(s')$. Then there exists $a\in N_GL$ such that $a$ is unipotent,
$as'=s'a\in\d$ and $s'a$ is isolated in $N_GL$.} 
\nl
We have ${}^\d\cz_L^0\sub Z_G(s'g_s)$ hence ${}^\d\cz_L^0\sub Z_G(s')$. Thus, 
$s\in Z_G({}^\d\cz_L^0)$. Since $Z_{G^0}({}^d\cz_L^0)=L$ (see 1.10) and 
$Z_{G^0}({}^\d\cz_L^0)$ is normal in $Z_G({}^\d\cz_L^0)$, we have $s'\in N_GL$. From
$s'g_s\in S_s$ we see that $g_s\in N_GL$. Since $g\in Z_G(s')$, we have 
$g_u\in Z_G(s')$. Now $Z_G(s')^0$ contains $g_s$ and ${}^\d\cz_L^0$; moreover, $g_s$
commutes with any element of ${}^\d\cz_L^0$ (since $s'g_s$ and $s'$ do). Hence we 
can find a maximal torus $T_1$ of $Z_G(s')^0$ such that $g_s\in T_1$,
${}^\d\cz_L^0\sub T_1$. Since $L=Z_{G^0}({}^\d\cz_L^0)$ we have $T_1\sub L$. Thus, 
$T_1$ is a torus in $Z_L(s')^0$. Since $g_s\in T_1$, we have $g_s\in Z_L(s')^0$. 
Since $s'g\in\bS U_P\sub N_GP$ we have $g_u=(s'g)_u\in N_GP$. Hence $g_u=ab$ where 
$a\in N_GL\cap N_GP,b\in U_P$ are uniquely determined and $a$ is unipotent. Now 
$s'g_s$ commutes with $g_u$ hence $(s'g_sag_s\i s'{}\i)(s'g_sbg_s\i s'{}\i)=ab$.
Since $s'g_s\in S_s\sub N_GL\cap N_GP$, we have
$a\i(s'g_s ag_s\i s'{}\i)=b(s'g_sbg_s\i s'{}\i)\i\in N_GL\cap U_P=\{1\}$. Hence 
$a\in Z_G(s'g_s)$. Since $Z_G(s'g_s)\sub Z_G(s')$ we have $a\in Z_G(s')$ and
$a\in Z_G(g_s)$. We have $s'g\in\bS U_P$ hence $s'g_sab\in\bS U_P$. Since 
$s'g_sa\in N_GL\cap N_GP,b\in U_P$, it follows that $s'g_sa\in\bS$. Let $g'=g_sa$. 
Then $g'\in N_GL$, $s'g'\in\bS$. We have $g'\in Z_{N_GL}(s')$. Since $a$ is 
unipotent we have $g'_s=g_s\in Z_L(s')^0$. Also, 
$Z_{N_GL}(s'g'_s)\sub Z_{N_GL}(s')$. Since $s'g'\in\bS$, we see that $s'g'$ is
isolated in $N_GL$. Applying (c) to $N_GL$ instead of $G$ we see that $s'a=as'$ is 
isolated in $N_GL$. Since $s'g_sa\in\bS$ we have $s'g_sa\in\d$. Since $g_s\in L$ we
have $s'a\in\d$. This proves (d).

(e) {\it Let $F$ be the image of $\{y\in\d;y\text{ isolated in }N_GL\}$ under
$y\m y_s$. Let $E$ be the set of all $g\in\d_1$ such that there exists $x\in G^0$ 
with $x\i sg_sx\in S_s,x\i sx\in F,x\i g_sx\notin{}^\d\cz_L^0$. Then $E$ is a closed
stable subset of $\d_1$.}
\nl
Let $\d'$ be the connected component of $N_GL$ such that $h\in\d\imp h_s\in\d'$.
From 2.7 and 1.22 we see that there exist finitely many semisimple $L$-conjugacy 
classes $C_0,C_1,\do,C_m$ in $\d'$ such that $F-S_s=\cup_{j=1}^m{}^\d\cz_L^0C_j$ and
$S_s={}^\d\cz_L^0C_0$. Applying an argument in 16.3(a) to ${}^\d\cz_L^0C_j$,
($j\in[0,m]$) instead of $S_s$ we see that there are only
finitely many orbits for the action $(h,l):x\m hxl\i$ of $Z_G(s)^0\T L$ on 
$\{x\in G^0;x\i sx\in{}^\d\cz_L^0C_j\}$. Hence for $j\in[0,m]$ we can find elements
$x_{ij}\in G^0$, $i\in[1,p_j],p_j<\iy$ such that 
$x_{ij}\i sx_{ij}\in{}^\d\cz_L^0C_j$ and $E=\cup_{j\in[0,m],i\in[1,p_j]}E_{ij}$ 
where
$$E_{ij}=\cup_{z\in Z_G(s)^0}\{g\in\d_1;z\i sg_sz\in x_{ij}{}^\d\cz_L^0C_0x_{ij}\i,
z\i g_sz\notin x_{ij}{}^\d\cz_L^0x_{ij}\i\}.$$
It is then enough to show that, for any $j,i$ as above, $E_{ij}$ is a closed stable
subset of $\d_1$. We set $a=x_{ij}$. We have $a\i sa\in\d'$ hence $a\i sa\in N_GL$.
Applying 1.27(a) with $N_GL,Z_{N_GL}(a\i sa),C_0$ instead of $H',H,\boc$ we see that
$C_0\cap Z_{N_GL}(a\i sa)=\cup_{r=1}^{m'}C'_r$ where $m'<\iy$ and $C'_r$ are 
semisimple $Z_L(a\i sa)^0$-conjugacy classes in $Z_{N_GL}(a\i sa)$. In the 
definition of $E_{ij}$ we have automatically

$z\i sg_sz\in a({}^\d\cz_L^0)(c_0\cap Z_{N_GL}(a\i sa))a\i$
\nl
(we use that ${}^\d\cz^0_L\sub Z_G(a\i sa)$ since $a\i sa\in\d'$). Hence
$E_{ij}=\cup_{r=1}^{m'}E_{ijr}$ where
$$E_{ijr}=\cup_{z\in Z_G(s)^0}\{g\in\d_1;z\i sg_sz\in\t aC'_ra\i,z\i g_sz\notin\t\}
$$
and $\t:=a({}^\d\cz_L^0)a\i$. It is then enough to show that, for any $r\in[1,m']$ 
as above, $E_{ijr}$ is a closed stable subset of $\d_1$. Let $f\in C'_r$. Since 
$aZ_L(a\i sa)^0a\i\sub Z_G(s)^0$, we have
$$E_{ijr}=\cup_{z\in Z_G(s)^0}\{g\in\d_1;z\i sg_sz\in\t afa\i,z\i g_sz\notin\t\}.$$
We may assume that $E_{ijr}$ is non-empty. Then $afa\i=sf_0$ where $f_0$ is a 
semisimple element of $Z_G(s)^0$. Now $f_0$ centralizes $\t$, a torus in 
$Z_G(s)^0$. Hence there exists a maximal torus $T_1$ of $Z_G(s)^0$ such that $T_1$ 
contains $\t$ and $f_0$. Applying (b) to $Z(s),\d_1,E_{ijr}$ instead of $G,D,\car$,
we see that it is enough to show that
$$\{h\in T_1;z\i shz\in sf_0\t,z\i hz\notin\t\text{ for some }z\in Z_G(s)^0\}$$
is closed in $T_1$. Since this is non-empty, we have $f_0\notin\t$ and the last 
variety becomes

$\{h\in T_1;z\i hz\in f_0\t\text{ for some }z\in Z_G(s)^0\}$
\nl
that is, $\cup_{w\in W_1}w\t f_0w\i$ where $W_1=N_{Z_G(s)^0}T_1/T_1$. This is 
clearly closed in $T_1$ since $W_1$ is finite and $\t f_0$ is closed in $T_1$. This
proves (e).

(f) {\it there exists an open stable subset $\cu_1$ of $\d_1$ such that $\cu_1$ 
contains any unipotent element in $\d_1$ and such that
$g\in\cu_1\imp Z_G(sg_s)\sub Z_G(s)$.}
\nl
We imbed $G$ into $\hG=GL_n(\kk)$ as a closed subgroup. Let
$\cu'_1=\{g\in\hG;Z_{\hG}(g_s)\sub Z_{\hG}(s)\}$. Let $\cu_1=s\i\cu'_1\cap\d_1$. 
Clearly, $\cu_1$ has the required properties.

\mpb

We can now prove the lemma. Let $\cu=\{g\in\cu_1;g\notin E\}$ with $\cu_1$ as in 
(f), $E$ as in (e). From (e),(f) we see 
that $\cu$ is an open stable subset of $\d_1$. If $g\in\d_1$ is unipotent, then 
$g\in\cu_1$ by (f) and $g\notin E$ (if we had $g\in E$ then there would exist 
$x\in G^0$ such that $1=x\i 1x\notin{}^\d\cz_L^0$, absurd). Thus $\cu$ contains any
unipotent element in $\d_1$. Assume now that $g\in\cu,x\in G^0,x\i sgx\in\bS U_P$ 
(with $P$ as in 16.3). We must show that $x\i sx\in S_sU_P$ and 
$x\i g_sx\in{}^\d\cz_L^0U_P$. Now any element in $\bS U_P$ is $U_P$-conjugate to an
element whose semisimple part is in $S_s$. (See the proof of 3.15 and 1.22(b).) 
Hence, replacing $x$ by $xv$ for some $v\in U_P$ we may assume that we have, in 
addition, $x\i sg_sx\in S_s$. Since $g\in\cu_1$, we have $Z_G(sg_s)\sub Z_G(s)$ 
hence $Z_G(xsg_sx\i)\sub Z_G(xsx\i)$. We apply (d) with $xsx\i,xgx\i$ instead of 
$s',g$. (We have $g\in Z_G(s)^0$ since $g\in\d_1$ and $\d_1$ contains unipotent 
elements. Hence $xgx\i\in Z_G(xsx\i)^0$.) We see that there exists $a\in N_GL$ such 
that $a$ is unipotent, $axsx\i=xsx\i a\in\d$ and $xsx\i a$ is isolated in $N_GL$. We
then have $xsx\i\in F$ ($F$ as in (e)). Since $x\i sg_sx\in S_s,x\i sx\in F$ and 
$g\notin E$ we must have $x\i g_sx\in{}^\d\cz_L^0$, by the definition of $E$. We 
have $x\i sx=(x\i sg_sx)(x\i g_sx)\i\in S_s{}^\d\cz_L^0=S_s$. This completes the 
proof of the lemma.

\subhead 16.5\endsubhead
Let $\cu$ be as in Lemma 16.4. Let $P$ be as in 16.3. Let $\ps:X@>>>\bY$ be as in 
3.14. We show:

(a) {\it the sets $X_{\cu,\co}=\{(g,xP)\in X;g\in s\cu,x\in\hco\}$ ($\co\in\G$) form
a finite partition of $X_\cu=\{(g,xP)\in X;g\in s\cu\}$ into open and closed 
subsets.} 
\nl
From 16.4(iii) we see that the second projection defines a morphism
$pr_2:X_\cu@>>>\{xP\in G^0/P;x\i sx\in S_sU_P\}$. We have 
$X_{\cu,\co}=pr_2\i(\hco/P)$ and it remains to use the fact that the subsets 
$\hco/P$ form a finite partition of $\{xP\in G^0/P;x\i sx\in S_sU_P\}$ into open and
closed subsets.

For any $x\in M$, $P_x:=xPx\i\cap Z_G(s)^0$ is a parabolic of $Z_G(s)^0$, see 16.3. 
Moreover, $L_x:=xLx\i\cap Z_G(s)^0$ is a Levi of $P_x$, since 
$x\i sx\in N_GL\cap N_GP$, see 1.12(a). Let 

$\boc_x=\{v\in Z_G(s);v\text{ unipotent, }x\i svx\in S\}$,

$\Si_x=x({}^\d\cz_L^0)x\i\boc_x$.
\nl
Let $\ce_x$ be the local system on $\Si_x$ obtained as the inverse image of $\ce$ 
under $\Si_x@>>>S,g\m x\i sgx$. The results of 15.1, 15.3, 15.4 are applicable to 
$Z_G(s),P_x,L_x,\Si_x,\ce_x$ instead of $G,P,L,\Si,\ce$ (see 16.2). Let

$\p_x:\tY'_x@>>>Y'_x,\tce_x,\ps_x:X'_x@>>>\bY'_x,K_x$
\nl
be obtained from $\p:\tY_{L,\Si}@>>>Y_{L,\Si},\tce,\ps:X@>>>\bY_{L,\Si},K$ in 15.1, 
15.3, 15.4, replacing $G,P,L,\Si,\ce$ by $Z_G(s),P_x,L_x,\Si_x,\ce_x$. 

For any $\co\in\G$ we choose a base point $x_\co\in\co$. We set

$P_\co=P_{x_\co},L_\co=L_{x_\co},\boc_\co=\boc_{x_\co},\Si_\co=\Si_{x_\co},
\ce_\co=\ce_{x_\co},\p_\co=\p_{x_\co}$,

$\tY'_\co=\tY'_{x_\co},Y'_\co=Y'_{x_\co},\tce_\co=\tce_{x_\co},\ps_\co=\ps_{x_\co},
X'_\co=X'_{x_\co},\bY'_\co=\bY'_{x_\co},K_\co=K_{x_\co}$,

$X'_{\cu,\co}=\{(h,zP_\co)\in\cu\T Z_G(s)^0/P_\co;z\i hz\in\baSi_\co U_{P_\co}\}\sub
X'_\co$.

\proclaim{Lemma 16.6} We have an isomorphism $X'_{\cu,\co}@>\si>>X_{\cu,\co}$,
$(h,zP_\co)\m(sh,zx_\co P)$.
\endproclaim 
We prove this only at the level of sets. First we show that our map is well defined.
Assume that $(h,zP_\co)\in X'_{\cu,\co}$. Then $z\i hz=\a\b\c$ with
$\a\in x_\co{}^\d\cz^0_Lx_\co\i,\b\in\bar\boc_\co,\c\in U_{P_\co}$ and
$$\align(zx_\co)\i shzx_\co&=x_\co\i s z\i hzx_\co=x_\co\i s\a\b\c x_\co\\&=(x_\co\i
\a x_\co)(x_\co\i s\b x_\co)(x_\co\i\c x_\co)\in{}^\d\cz_L^0\bS U_P=\bS U_P,
\endalign$$
$sh\in s\cu,zx_\co\in\co$, hence $(sh,zx_\co P)\in X_{\cu,\co}$. (We use that 
$x_\co\i s\bar\boc_\co x_\co\sub\bS$, $x_\co\i U_{P_\co}x_\co\sub U_P$.)

We show that our map is injective. Assume that 
$(h,zP_\co),(h',z'P_\co)\in\cu\T Z_G(s)^0/P_\co$ satisfy 
$(sh,zx_\co P)=(sh',z'x_\co P)$. Then clearly $h=h'$ and 
$z\i z'\in x_\co Px_\co\i\cap Z_G(s)^0=P_\co$ hence $zP_\co=z'P_\co$.

We show that our map is surjective. Let $(g,xP)\in X_{\cu,\co}$. We have $g=sh$,
$xP=zx_\co P$ where 

$h\in\cu,z\in Z_G(s)^0$, $x_\co\i z\i shzx_\co\in\bS U_P$. 
\nl
We have $z\i hz=a'c$ where $a'\in s\i x_\co\bS x_\co\i,c\in U_{x_\co Px_\co\i}$. 
Since $z\i hz\in Z_G(s)$ we must have $a'\in Z_G(s),c\in Z_G(s)$. Now 

$c\in U_{x_\co Px_\co\i}\cap Z_G(s)=U_{x_\co Px_\co\i}\cap Z_G(s)^0=U_{P_\co}$,
\nl
by 1.11. Thus $z\i hz\in a'U_{x_\co Px_\co\i}$ where 
$a'\in(s\i x_\co\bS x_\co\i)\sub N_G(x_\co Lx_\co\i)$. Hence 

$(z\i hz)_s\in a'_s U_{P_\co}$ with $a'_s\in N_G(x_\co Lx_\co\i)$.
\nl
Using 16.4(iii) we have $x_\co\i z\i h_szx_\co\in{}^\d\cz_L^0U_P$. Thus, 
$(z\i hz)_s=z\i h_sz\in a(x_\co U_Px_\co\i)$ with $a\in x_\co{}^\d\cz_L^0x_\co\i$.
Since 

$a\i a'_s\in N_G(x_\co Lx_\co\i)\cap U_{x_\co Px_\co\i}=\{1\}$
\nl
we have $a'_s=a\in x_\co{}^\d\cz_L^0x_\co\i$. Let $b=a'_u$. It remains to show that
$b\in\bar\boc_\co$. Since $a'\in Z_G(s)$, we have $b\in Z_G(s)$. From
$sa'\in x_\co\bS x_\co\i$ we deduce 

$sa'_u\in x_\co\bS x_\co\i a'_s{}\i\sub x_\co\bS x_\co\i x_\co{}^\d\cz_L^0x_\co\i
=x_\co\bS x_\co\i$. 
\nl
It follows that $a'_u\in\bar\boc_\co$, as required. The lemma is proved.

\subhead 16.7\endsubhead
Let $\tY=\tY_{L,S},Y=Y_{L,S}$ and let $\p:\tY@>>>Y$ be as in 3.13. Let 
$\tY_\cu=\{(g,xL)\in\tY;g\in s\cu\}$. For any $\co\in\G$ we set

$\tY_{\cu,\co}=\{(g,xL)\in\tY;g\in s\cu,x\in\co\},\quad
Y_{\cu,\co}=\p(\tY_{\cu,\co})$,

$\tY'_{\cu,\co}=\{(h,zL_\co)\in\tY'_\co;h\in\cu\}$.

\proclaim{Lemma 16.8}(a) The map $(g,xL)\m(g,xP)$ is an isomorphism of $\tY_\cu$ 
onto the open subset $\ps\i(Y\cap s\cu)$ of $X_\cu$.

(b) The subsets $\tY_{\cu,\co},(\co\in\G)$ form a finite partition of $\tY_\cu$ into
open and closed subsets.

(c) The map $\ps:X_\cu@>>>\bY\cap s\cu$ is proper, surjective, and $Y\cap s\cu$ is
open in $\bY\cap s\cu$.

(d) The map $(h,zL_\co)\m(h,zP_\co)$ is an isomorphism of $\tY'_{\cu,\co}$ onto 
the open subset $\ps_\co\i(Y'_\co\cap\cu)$ of $X'_{\cu,\co}$. The map $\ps_\co$ 
restricts to a proper map of $X'_{\cu,\co}$ onto $\bY'_\co\cap\cu$ and 
$Y'_\co\cap\cu$ is open, dense in $\bY'_\co\cap\cu$.
\endproclaim
We prove (a). By 5.5, the same formula gives an isomorphism $\tY@>\si>>\ps\i(Y)$.
Hence the map in (a) is an isomorphism onto $\ps\i(Y\cap s\cu)$. It remains to show
that $\ps\i(Y\cap s\cu)$ is open in $X_\cu$. Since 
$\ps\i(Y\cap s\cu)=\ps\i(Y)\cap X_\cu$ it is enough to show that $\ps\i(Y)$ is open
in $X$. This follows from the fact that $Y$ is open in
$\bY=\ps(X)$.

We prove (b). The map in (a) identifies $\tY_\cu$ with an open subset of $X_\cu$ and
$\tY_{\cu,\co}$ with $\tY_\cu\cap X_{\cu,\co}$. This together with 16.5(a) yields 
(b).

We prove (c). We have $X_\cu=\ps\i(\bY\cap s\cu)$ hence the first assertion of (b)
follows by change of base from the fact that $\ps$ is proper. Since $Y$ is open in 
$\bY$, we see that $Y\cap(\bY\cap s\cu)$ is open in $\bY\cap s\cu$. Hence 
$Y\cap s\cu$ is open in $\bY\cap s\cu$. 

We prove (d). From 15.3(a) we see that $(h,zL_\co)\m(h,zP_\co)$ gives an isomorphism
of $\tY'_\co$ onto the open subset $\ps_\co\i(Y'_\co)$ of $X'_\co$. Hence the same 
formula gives an isomorphism of $\tY'_{\co,\cu}$ onto the open subset 
$\ps_\co\i(Y'_\co\cap\cu)$ of $X'_{\cu,\co}$. The map $\ps_\co:X'_\co@>>>\bY'_\co$
is proper, surjective. Since $\ps_\co:X'_{\cu,\co}@>>>\bY'_\co\cap\cu$ is obtained 
from the previous map by change of base, it is also proper, surjective. Since 
$Y'_\co$ is open in $\bY'_\co$, $Y'_\co\cap\cu$ is open in $\bY'_\co\cap\cu$. We 
show it is dense. We have $\bY'_\co=\cup_F(F\cap\cu)$ 
where $F$ runs over the irreducible components of $\bY'_\co$. It is enough to show 
that $F\cap Y'_\co\cap\cu$ is dense in $F\cap\cu$ for any $F$. We may assume that 
$F\cap\cu\ne\em$. Since $F\cap\cu$ is open, non-empty in the irreducible variety 
$F$, it is also dense in $F$. Since $Y'_\co$ is open dense in $\bY'_\co$ (see 15.1),
we see that $F\cap Y'_\co$ is open dense in $F$. Since $F\cap\cu$, $F\cap Y'_\co$ 
are open dense in $F$, their intersection $F\cap Y'_\co\cap\cu$ is open dense in 
$F\cap\cu$. The lemma is proved.

\proclaim{Lemma 16.9}(a) The map $\ps:X_{\cu,\co}@>>>\bY\cap s\cu$ is proper, with 
image equal to $s(\bY'_\co\cap\cu)$.

(b) Let us identify $\tY_{\cu,\co}$ with a subset of $X$ via the imbedding 
$\tY\sub X$ in 5.5. Then $\tY_{\cu,\co}=\ps\i(Y_{\cu,\co})\cap X_{\cu,\co}$.

(c) $Y_{\cu,\co}$ is open in $s(\bY'_\co\cap\cu)$.

(d) $Y_{\cu,\co}$ is open and closed in $Y\cap s\cu$. We have
$\cup_{\co\in\G}Y_{\cu,\co}=Y\cap s\cu$. For $\co,\co'\in\G$,
$Y_{\cu,\co},Y_{\cu,\co'}$ coincide if $\co,\co'$ are in the same $\cw_S$-orbit in
$\G$ and are disjoint, otherwise.

(e) For any $\co\in\G$, $\tY_{\cu,\co}$ is a dense subset of $X_{\cu,\co}$ (see 
(b)).

(f) For any $\co\in\G$, $Y_{\cu,\co}$ is an open dense subset of 
$s(\bY'_\co\cap\cu)$.
\endproclaim
We prove (a). The fact that $\ps:X_{\cu,\co}@>>>\bY\cap s\cu$ is proper follows from
16.8(c) since $X_{\cu,\co}$ is closed in $X_\cu$ (see 16.5(a)). The statement about 
its image can be reduced using the isomorphism in 16.6 to a statement in 16.8(d).

We prove (b). We must show that $\tY_{\cu,\co}$ is a union of fibres of
$\ps:X_{\cu,\co}@>>>s(\bY'_\co\cap\cu)$. This is clear from the definitions.

We prove (c). From the proof of (b) we see that $X_{\cu,\co}-\tY_{\cu,\co}$ is also
a union of fibres of $\ps:X_{\cu,\co}@>>>s(\bY'_\co\cap\cu)$. Hence its image under
the proper surjective map $\ps:X_{\cu,\co}@>>>s(\bY'_\co\cap\cu)$ is a closed subset
of $s(\bY'_\co\cap\cu)$ complementary to the image $Y_{\cu,\co}$ of $\tY_{\cu,\co}$.
This proves (c).

We prove (d). The map $\p:\tY_\cu@>>>Y\cap s\cu$ is proper, surjective, since it is 
obtained by change of base from the proper surjective map $\p:\tY@>>>Y$. Since 
$\tY_{\cu,\co}$ is closed in $\tY_\cu$ (see 16.8(b)), it follows that 
$Y_{\cu,\co}=\p(\tY_{\cu,\co})$ is closed in $Y\cap s\cu$. Since 
$\cup_\co\tY_{\cu,\co}=\tY_\cu$ (see 16.8(b)), it follows that
$\cup_\co Y_{\cu,\co}=Y\cap s\cu$.

Assume that $Y_{\cu,\co},Y_{\cu,\co'}$ are not disjoint; let 
$g\in Y_{\cu,\co}\cap Y_{\cu,\co'}$. Then there exist $x\in\co,x'\in\co'$ such that
$(g,xL)\in\tY,(g,x'L)\in\tY$. Using 3.13(a), we see that there exists 
$n\in N_{G^0}L$ such that $nSn\i=S$ and $x'=xn\i$. Then $(g_1,x_1L)\m(g_1,x_1n\i L)$
is a bijection $\tY_{\cu,\co}@>>>\tY_{\cu,\co'}$ hence $Y_{\cu,\co}=Y_{\cu,\co'}$. 
The same argument shows that, if $\co,\co'$ are in the same $\cw_S$-orbit then 
$Y_{\cu,\co}=Y_{\cu,\co'}$. We see that the complement of $Y_{\cu,\co}$ in 
$Y\cap s\cu$ is the union of the closed subsets $Y_{\cu,\co'}$ (with $\co'$ not in 
the $\cw_S$-orbit of $\co$); hence $Y_{\cu,\co}$ is open in $Y\cap s\cu$. This 
proves (d).

We prove (e). We have $\boc_\co=\sqc_{j\in J}\g_j$ where $J$ is finite and $\g_j$ 
are (unipotent) $L_\co$-conjugacy classes in $N_{Z_G(s)}L_\co$. We have 
$X'_\co=\cup_{j\in J}X'_j,\tY'_\co=\sqc_{j\in J}\tY'_j$ where
$$X'_j=\{(h,zP_\co)\in Z_G(s)\T Z_G(s)^0;
z\i hz\in x_\co{}^\d\cz_L^0x_\co\i\bag_jU_{P_\co}\},$$
$$\align\tY'_j=\{&(h,zL_\co)\in Z_G(s)\T Z_G(s)^0;\\&
z\i hz=ab,a\in x_\co{}^\d\cz_L^0x_\co\i,
b\in\g_j,Z_{Z_G(s)}(a)^0\sub L_\co\}.\endalign$$
Let $X'_{\cu,j}=\{(h,zP_\co)\in X'_j;h\in\cu\},\tY'_{\cu,j}=
\{(h,zL_\co)\in\tY'_j;h\in\cu\}$. 

We have $X'_{\cu,\co}=\cup_{j\in J}X'_{\cu,j}$,
$\tY'_{\cu,\co}=\cup_{j\in J}\tY'_{\cu,j}$.

Since $X'_j$ is irreducible and $\cu$ is open in $Z_G(s)$ we see that $X'_{\cu,j}$
is open in $X'_j$ so that $X'_{\cu,j}$ is either empty or irreducible. If 
$X'_{\cu,j}\ne\em$ then we can find
$a\in x_\co{}^\d\cz_L^0x_\co\i,b\in\bag_j,c\in U_{P_\co}$ with $abc\in\cu$. Since
$a$ and $c$ are contained in $Z_G(s)^0$, we see that $b$ is contained in $\d_1$, the
connected component of $Z_G(s)$ that contains $\cu$. Hence $\g_j\sub\d_1$. Let us
identify $X'_{\cu,j}$ with a subset $X_{\cu,j}$ of $X_{\cu,\co}$ via the isomorphism
$X'_{\cu,\co}@>>>X_{\cu,\co}$ in 16.6. Since 
$X_{\cu,\co}=\cup_{j\in J;\g_j\sub\d_1}X_{\cu,j}$, we see that it is enough to show 
that, if $\g_j\sub\d_1$, then

(g) $\tY_{\cu,\co}\cap\tY'_{\cu,j}$ is dense in $\tY'_{\cu,j}$,

(h) $\tY'_{\cu,j}$ is dense in $X'_{\cu,j}=X_{\cu,j}$.
\nl
(We regard $\tY'_{\cu,j}$ as a subspace of $X'_{\cu,j}$ via the imbedding
$\tY'_{\cu,\co}@>>>X'_{\cu,\co}$ in 16.8(d) hence as a subspace of $X_{\cu,j}$.) 

We prove (g). This is equivalent to the following statement:
$$\align\{(h,zL_\co);&h\in\cu,z\in Z_G(s)^0,z\i hz=ab,
a\in x_\co{}^\d\cz_L^0x_\co\i,\\&b\in\g_j,(Z(s)\cap Z(a))^0\sub L_\co,
Z_G(sa)^0\sub x_\co Lx_\co\i\}\endalign$$
is dense in
$$\align\{(h,zL_\co);&h\in\cu,z\in Z_G(s)^0,z\i hz=ab, 
a\in x_\co{}^\d\cz_L^0x_\co\i,\\&b\in\g_j, (Z(s)\cap Z(a))^0\sub L_\co\}.
\endalign$$
Since the condition $Z_G(sa)^0\sub x_\co Lx_\co\i$ implies that
$(Z(s)\cap Z(a))^0\sub L_\co$, we see that it is enough to show that for any 
$b\in\g_j$,

$\{a\in x_\co{}^\d\cz_L^0x_\co\i;Z_G(sa)^0\sub x_\co Lx_\co\i,ab\in\cu\}$ is dense
in $x_\co{}^\d\cz_L^0x_\co\i$
\nl
or that

$\{a\in x_\co{}^\d\cz_L^0x_\co\i;Z_G(sa)^0\sub x_\co Lx_\co\i\}
\cap(x_\co{}^\d\cz_L^0x_\co\i\cap\cu b\i)$ is dense in $x_\co{}^\d\cz_L^0x_\co\i$.
\nl
Since $b\in\cu$, $x_\co{}^\d\cz_L^0x_\co\i\cap\cu b\i$ is an open subset of the 
torus $x_\co{}^\d\cz_L^0x_\co\i$ containing the unit element; hence it is an open 
dense subset of $x_\co{}^\d\cz_L^0x_\co\i$. On the other hand, 
$\{a\in x_\co{}^\d\cz_L^0x_\co\i;Z_G(sa)^0\sub x_\co Lx_\co\i\}$ is an open dense
subset of $x_\co{}^\d\cz_L^0x_\co\i$ by 3.10(a). Since the intersection of two open 
dense subsets of a torus is dense in that torus, (g) is proved.

We prove (h). Since $\tY'_{\cu,j}$ is open in $\tY'_j$ which is open in $X'_j$, we 
see that $\tY'_{\cu,j}$ is open in $X'_j$. Since $X'_j$ is irreducible, to prove (h)
it suffices to show that $\tY'_{\cu,j}\ne\em$. But this is contained in the proof of
(g). This proves (h) hence (e).

We prove (f). The openness follows from (c). From (e) we see that 
$\ps(\tY_{\cu,\co})$ is a dense subset of $\ps(X_{\cu,\co})$ hence $Y_{\cu,\co}$ is
a dense subset of $s(\bY'_\co\cap\cu)$ (see 16.9(a)). The lemma is proved.

\subhead 16.10\endsubhead
For a $\cw_S$-orbit $Z$ in $\G$ we set $Y_{\cu,Z}=Y_{\cu,\co}$, 
$\bY'_{Z,\cu}=\bY'_\co\cap\cu$ where $\co\in\G$. This is well defined, by
16.9(d),(f). For $\co\in\G$, 

(a) $s(Y'_\co\cap\cu)\cap Y_{\cu,\co}$ is open in $Y_{\cu,\co})$. 
\nl
(Since $Y_{\cu,\co}\sub s(\bY'_\co\cap\cu)$, it suffices to show that 
$s(Y'_\co\cap\cu)$ is open in $s(\bY'_\co\cap\cu)$. This follows from 16.8(d).) 
Hence, for any $\cw_S$-orbit $Z$ in $\co$,
$$\cv_Z:=\cap_{\co\in Z}(s(Y'_\co\cap\cu)\cap Y_{\cu,Z})$$
is an open subset of $Y_{\cu,Z}$. It also follows that $\cv_Z$ is open in 
$Y\cap s\cu$. Let $\cv=\cup_Z\cv_Z$ where $Z$ runs over the $\cw_S$-orbits in $\G$.
We show:

(b) {\it $\cv$ is an open smooth subset of $Y\cap s\cu$ of pure dimension
$\dim Z_G(s)^0-\dim L+\dim S$ and the $\cv_Z$ form a finite partition of $\cv$ into 
open and closed subsets.}
\nl
Since $\cv_Z$ is open in $Y\cap s\cu$, the union $\cv=\cup_Z\cv_Z$ is open in 
$Y\cap s\cu$ and $\cv_Z$ is open in $\cv$ for all $Z$. For $Z\ne Z'$, the sets
$\cv_Z,\cv_{Z'}$ are disjoint since they are contained in $Y_{\cu,Z},Y_{\cu,Z'}$ 
which are disjoint (see 16.9(d)). Hence the sets $\cv_Z$ are also closed in $\cv$.
For any $Z$ and any $\co\in Z$, $\cv_Z$ is open in $s(Y'_\co\cap\cu)$. (It is 
enough to show that $Y_{\cu,Z}\cap\cap_{\co'\in Z}s(Y'_{\co'}\cap\cu)$ is open in 
$s(\bY'_\co\cap\cu)$. This follows from the fact that $Y_{\cu,\co}$ is open in 
$s(\bY'_\co\cap\cu)$, see 16.9(c), and that
$Y_{\cu,Z}\cap\cap_{\co'\in Z}s(Y'_{\co'}\cap\cu)$ is open in $Y_{\cu,Z}$, see 
(a).) 
Since $s(Y'_\co\cap\cu)$ is isomorphic to an open set in $Y'_\co$, it follows that
$\cv_Z$ is isomorphic to an open set in $Y'_\co$. Since $Y'_\co$ is smooth of pure
dimension $\dim Z_G(s)^0-\dim L_\co+\dim {}^\d\cz_L^0+\dim\boc_\co$ (see 15.1 and 
3.13(b)), it follows that $\cv_Z$ is smooth of pure dimension 
$\dim Z_G(s)^0-\dim L_\co+\dim{}^\d\cz_L^0+\dim\boc_\co$ where $\co\in Z$. Now 
$L_\co$ is isomorphic to the connected centralizer in $L$ of an element in $S_s$ 
hence $\dim L_\co=\dim L-(\dim S_s-\dim{}^\d\cz^0_L)$; moreover,
$\dim\boc_\co=\dim S-\dim S_s$. Hence $\cv_Z$ is smooth of pure dimension 
$\dim Z_G(s)^0-\dim L+\dim S$ which is independent of $Z$. The same is then true for
$\cv$. This proves (b).

We show:

(c) {\it $\cv$ is open dense in $\bY\cap s\cu$.}
\nl
By (b), $\cv$ is open in $Y\cap s\cu$ and by 16.8(c), $Y\cap s\cu$ is open in 
$\bY\cap s\cu$. Hence $\cv$ is open in $\bY\cap s\cu$. We prove that it is also 
dense. We have
$$\bY\cap s\cu=\ps(X_\cu)=\cup_\co\ps(X_{\cu,\co})=s\cup_\co\ps_\co(X'_{\cu,\co})=
s\cup_\co(\bY'_\co\cap\cu)=\cup_Z s\bY'_{Z,\cu}$$
where $Z$ runs over the $\cw_S$-orbits in $\G$ and $\bY'_{Z,\cu}$ is as in 16.10. 
(We have used 16.8(c), 16.5(a), 16.6, 16.8(d).) Since $\cv=\cup_Z\cv_Z$, it is 
enough to show that, for any $Z$, $\cv_Z$ is dense in $s\bY'_{Z,\cu}$. By 16.8(d) 
and 16.9(f), for any $\co\in\cz$, $s(Y'_\co\cap\cu)$ is open dense in 
$s\bY'_{Z,\cu}$ and $Y_{\cu,\co}$ is open dense in $s\bY'_{Z,\cu}$. Hence 
$s(Y'_\co\cap\cu)\cap Y_{\cu,\co}$ is open dense in $s\bY'_{Z,\cu}$. Hence 
$\cv_Z=\cap_{\co\in Z}(s(Y'_\co\cap\cu)\cap Y_{\cu,\co})$ is open dense in 
$s\bY'_{Z,\cu}$. This proves (c).

\proclaim{Lemma 16.11} Let ${}^0\tY=\{(g,xL)\in\tY;g\in\cv\}$. For $\co\in\G$ let

${}^0\tY'_\co=\{(h,zL_\co)\in\tY'_\co;h\in s\i\cv_Z\}$ 
\nl
where $\co\in Z$. We have a well defined isomorphism
$\a:\sqc_{\co\in\G}{}^0\tY'_\co@>\si>>{}^0\tY,(h,zL_\co)\m(sg,zx_\co L)$.
\endproclaim
For $\co\in\G$ let ${}^0\tY_\co=\{(g,xL)\in\tY;g\in\cv,x\in\co\}$. Since 
$\cv\sub s\cu$, the subsets ${}^0\tY_\co$ form a partition of ${}^0\tY$ into open
and closed subsets (see 16.8(b)). It is enough to show that for any $\co$, we 
have a well defined isomorphism 
${}^0\tY'_\co@>\si>>{}^0\tY_\co,(h,zL_\co)\m(sg,zx_\co L)$. The imbedding 
$\tY_\cu@>>>X_\cu$ (see 16.8(a)) identifies ${}^0\tY_\co$ with an open subset of 
$X_{\cu,\co}$; the imbedding $\tY'_{\cu,\co}@>>>X'_{\cu,\co}$ identifies 
${}^0\tY'_\co$ with an open subset of $X'_{\cu,\co}$. It is enough to show that the
isomorphism $X'_{\cu,\co}@>\si>>X_{\cu,\co}$ in 16.6 carries the subspace 
$\tY'_{\cu,\co}$ onto the subspace ${}^0\tY_\co$. Thus, it is enough to show that,
for $(h,zP_\co)\in X'_\co$, the following two conditions are equivalent:

(i) $h\in s\i\cv_Z, h\in Y'_\co$;

(ii) $sh\in\cv, sh\in Y_{\cu,\co}$.
\nl
Both (i) and (ii) are equivalent to the condition $sh\in\cv_Z$. The lemma is proved.

\subhead 16.12\endsubhead
By 16.11 we have a commutative diagram
$$\CD
\sqc_{\co\in\G}{}^0\tY'_\co@>\a>>{}^0\tY\\
@VVV            @VVV \\
s\i\cv @>\e>>\cv
\endCD$$
where $\e(h)=sh$ and the vertical maps are given by the first projection. Hence we
have a canonical isomorphism 
$$\op_{\co\in\G}(\p_{\co !}\tce_\co)|_{s\i\cv}@>\si>>\e^*((\p_!\tce)|_\cv)$$
of local systems on $s\i\cv$. (Here $(\p_{\co !}\tce_\co)|_{s\i\cv}$ is by 
definition the restriction of $\p_{\co !}\tce_\co$ to $s\i\cv_Z$ where $\co\in Z$ 
and is zero on $s\i\cv_{Z'}$ for $Z'\ne Z$.) This can be also regarded as an 
isomorphism
$$\op_{\co\in\G}IC(\bY'_\co,\p_{\co !}\tce_\co)|_{s\i\cv}@>\si>>
\e^*(IC(\bY,\p_!\tce)|_\cv).\tag a$$
Assume that we can show that the isomorphism (a) is the restriction to $s\i\cv$ of
an isomorphism
$$\op_{\co\in\G}IC(\bY'_\co,\p_{\co !}\tce_\co)|_{\bY'_\co\cap\cu}@>\si>>
\e^*(IC(\bY,\p_!\tce)|_{\bY\cap\cu}).\tag b$$
(Here $\e$ is regarded as an isomorphism $s\i\bY\cap\cu@>>>\bY\cap s\cu,g\m sg$;
moreover, $IC(\bY'_\co,\p_{\co !}\tce_\co)|_{\bY'_\co\cap\cu}$ is regarded as a 
complex on $s\i\bY\cap\cu$, zero outside $\bY'_\co\cap\cu$). An isomorphism (b) 
extending (a) is unique, if it exists. (This follows from the fact that the left    
hand side of (b) is the intersection cohomology complex of $s\i\bY\cap\cu$ with 
coefficients in a local system on the open dense smooth subvariety $s\i\cv$ of pure
dimension,  namely $\op_{\co\in\G}(\p_{\co !}\tce_\co)|_{s\i\cv}$.) The isomorphism
(b) gives rise for any $i$ to an isomorphism of stalks
$$\op_{\co\in\G}\ch^i_uIC(\bY'_\co,\p_{\co !}\tce_\co)@>\si>>
\ch^i_{su}IC(\bY,\p_!\tce).\tag c$$
We have $u\in s\i\bY\cap\cu$. Indeed, $su\in\bY$ and $u\in\cu$ since $u\in\d_1$ and
$\cu$ contains any unipotent element in $\d_1$.

We now show the existence of the isomorphism (b). To do this we will use $P$ in
16.3. (However, the isomorphism we construct will be independent of the choice of
$P$ in view of its uniqueness.) Using Lemma 5.7 and 15.4(a) we find isomorphisms
$$IC(\bY,\p_!\tce)|_{\bY\cap s\cu}@>\si>>\ps_!(K|_{X_\cu}) 
\text{ in }\cd(\bY\cap s\cu),\tag d$$
$$IC(\bY'_\co,\p_{\co !}\tce_\co)|_{\bY'_\co\cap\cu}@>\si>>
\ps_{\co !}(K_\co|_{X'_{\cu,\co}})\text{ in }\cd(\bY'_\co\cap\cu),\tag e$$
($K$ as in 5.7.) From 16.5(a) and 16.6 we get an isomorphism
$$\op_{\co\in\G}\ps_{\co !}(K_\co|_{X'_{\cu,\co}})@>\si>>\e^*(\ps_!(K|_{X_\cu}))
\text{ in }\cd(s\i\bY\cap\cu).\tag f$$
(We regard $\ps_{\co !}(K_\co|_{X'_{\cu,\co}})$ as a complex on $s\i\bY\cap\cu$
equal to $0$ outside $\bY'_\co\cap\cu$.) Combining the isomorphisms (d),(e),(f) we 
obtain an isomorphism as in (b).

\subhead 16.13\endsubhead
In the remainder of this section we assume that $\kk$ is an algebraic closure of a 
finite field $\FF_q$ and that $G$ has a fixed $\FF_q$-structure with Frobenius map 
$F:G@>>>G$. Assume that $F(L)=L,F(S)=S,\ce\in\cs(S)$ and that we are given an 
isomorphism $\ph_0:F^*\ce@>\si>>\ce$ of local systems on $S$. Then $F(\bY)=\bY$ and
$\ph_0$ induces an isomorphism $\ph:F^*\fK@>\si>>\fK$ where $\fK=IC(\bY,\p_!\tce)$.
Assume that $s,u$ in 16.3 satisfy $F(s)=s,F(u)=u$. For any $x\in(G^0)^F$ such that 
$x\i sx\in S_s$ let $L_x,\boc_x,\Si_x,\ce_x$ be as in 16.5. Let 
$\cf_x=\ce_x|_{\boc_x}$, a local system on $\boc_x$. Now $Z_G(s),L_x,\boc_x,\Si_x$
are defined over $\FF_q$ and $\ph_0:F^*\ce@>\si>>\ce$ induces an isomorphism 
$\ph'_x:F^*\ce_x@>\si>>\ce_x$ and an isomorphism $\ph_x:F^*\cf_x@>\si>>\cf_x$. Then
$L_x,Z_G(s),\boc_x,\Si_x,\cf_x,\ph_x$ are like $L,G,\Si^\o,\Si,\cf,\ph_1$ in 15.12 
hence the generalized Green function

$Q_{L_x,Z_G(s),\boc_x,\cf_x,\ph_x}:\{\text{unipotent elements in } Z_G(s)^F\}@>>>
\bbq$
\nl
is defined as in 15.12(c). We have the following result.

\proclaim{Theorem 16.14} $$\c_{\fK,\ph}(su)=\sum_{x\in(G^0)^F;x\i sx\in S_s}
\fra{|L_x^F|}{|Z_G^0(s)^F||L^F|}Q_{L_x,Z_G(s),\boc_x,\cf_x,\ph_x}(u).$$
\endproclaim
We can choose the base points $x_\co$ in $\co$ (see 16.5) in such a way that
$F(x_\co)=x_{F(\co)}$ for any $\co\in\G$. (We use the fact that, if $\co$ is
$F^n$-stable, then $\co^{F^n}\ne\em$. This follows from the fact that $\co$ is a
homogeneous space under a connected group.) Now the sum over $x$ in the theorem can
be broken into sums over $x\in\co^F$ for various $\co\in\G$ with $F(\co)=\co$. The 
sum over $x\in\co^F$ is equal to $Q_{L_\co,Z_G(s),\boc_\co,\cf_\co,\ph_\co}(u)$ 
(notation of 16.5) since all terms of the sum are equal and the number of terms in
the sum is $|\co^F|=\fra{|Z_G^0(s)^F||L^F|}{|L_\co^F|}$. Thus the right hand side of
the equality in the theorem is 

$\sum_{\co\in\G;F(\co)=\co}Q_{L_\co,Z_G(s),\boc_\co,\cf_\co,\ph_\co}(u)$
\nl
or equivalently,

$\sum_{\co\in\G;F(\co)=\co}\c_{\fK_\co,\ti\ph}(u)$
\nl
where $\fK_\co$ is $IC(\bY'_\co,\p_{\co !}\tce_\co)$ (notation of 16.5) extended by
$0$ on $\bY-\bY'_\co$ and $\ti\ph:F^*\fK_\co@>\si>>\fK_\co$ is induced by
$\ph'_x:F^*\ce_x@>\si>>\ce_x$. (At this point, it is important that
$Q_{L_\co,Z_G(s),\boc_\co,\cf_\co,\ph_\co}(u)$ can be computed in terms of a not 
necessarily "trivial choice", see 15.12(d).) We see that it is enough to prove that

(a) $\c_{\fK,\ph}(su)=\sum_{\co\in\G;F(\co)=\co}\c_{\fK_\co,\ti\ph}(u)$.
\nl
We shall now make the choice of $\cu$ in Lemma 16.4 more precise. Namely we will 
choose it so that, in addition, it satisfies $F(\cu)=\cu$. (In the proof of 16.4(f)
we choose the imbedding $G\sub GL_n(\kk)$ so that it is defined over $\FF_q$. Then 
$\cu_1$ defined in that proof is automatically $F$-stable hence 
$\cu=\{g\in\cu_1;g\notin E\}$, see the proof of 16.4, is again $F$-stable.) With 
this choice of $\cu$, the isomorphism 16.12(c) commutes with the natural Frobenius 
maps on its two sides. This gives rise to the equality of (alternating sums of) 
traces (a). The theorem is proved.

\enddocument